

\input amstex

\documentstyle{amsppt}

\loadbold

\magnification=\magstep1

\pageheight{9.0truein}
\pagewidth{6.5truein}

\NoBlackBoxes

\def\ann{\operatorname{ann}}
\def\spec{\operatorname{Spec}}
\def\SPEC{\operatorname{SPEC}}
\def\Mod{\operatorname{Mod}}
\def\rank{\operatorname{rank}}
\def\r{\bold{r}}
\def\c{\bold{c}}

\def\closed{\operatorname{Closed}}

\def\ArtSch{1}
\def\ArtZha{2}
\def\Bra{3}
\def\Dix{4}
\def\Gol{5}
\def\GooWar{6}
\def\Letone{7}
\def\Lettwo{8}
\def\McCRob{9}
\def\NdiVOy{10}
\def\Rosone{11}
\def\Rostwo{12}
\def\Smione{13}
\def\Smitwo{14}
\def\SmiZha{15}
\def\Ste{16}
\def\VdB{17}
\def\VOy{18}
\def\VOyVerone{19}
\def\VOyVertwo{20}
\def\War{21}

\topmatter

\title On Continuous and Adjoint Morphisms Between Noncommutative Prime
Spectra \endtitle

\rightheadtext{Noncommutative Spectra}

\author Edward S. Letzter \endauthor

\abstract We study topological properties of the correspondence of
prime spectra associated to a noncommutative ring homomorphism $R
\rightarrow S$. Our main result provides criteria for the adjointness
of certain functors between the categories of Zariski closed subsets
of $\spec R$ and $\spec S$; these functors arise naturally from
restriction and extension of scalars. When $R$ and $S$ are left
noetherian, adjointness occurs only for centralizing and ``nearly
centralizing'' homomorphisms. \endabstract



\address Department of Mathematics, Temple University, Philadelphia,
PA 19122 \endaddress

\email letzter\@math.temple.edu\endemail

\thanks The author's research was supported in part by NSF grants
DMS-9970413 and DMS-0196236. \endthanks

\endtopmatter

\document

\head 1. Introduction \endhead 

One of the most elementary and well-known properties of noncommutative rings
is the non-functoriality of their prime spectra: There is apparently no
natural way of assigning, to an arbitrary ring homomorphism $R \rightarrow S$,
a function from the prime spectrum of $S$ into the prime spectrum of $R$.
Nevertheless, there is an extensive and deep literature presenting -- among
many other things -- topological and geometric contexts for both
noncommutative ring homomorphisms and their generalizations to certain
functors between module-like categories.  These contexts appear, for example,
in the earlier \cite{\ArtSch; \Gol; \VOy; \VOyVerone; \VOyVertwo} and the more
recent \cite{\ArtZha, \NdiVOy; \Rosone; \Rostwo; \Smione; \Smitwo;
\SmiZha}. In the present paper we continue a discussion begun in
\cite{\ArtSch, \S 4}. We focus on topological properties of the
correspondences of prime spectra associated to arbitrary homomorphisms
involving left noetherian rings or affine PI algebras.

\subhead 1.1 \endsubhead To fix notation, equip the set $\spec R$ of
prime ideals of a (not necessarily commutative) ring $R$ with the
Zariski topology, by declaring the closed subsets to be those of the
form
$$V_R(X) = \{ P \in \spec R : P \supseteq X \},$$
for $X \subseteq R$. Our specific intent in this paper is to carefully
examine noncommutative generalizations of the following two trivially
true but fundamentally important facts: If $f\colon R \rightarrow S$
is a commutative ring homomorphism then (1) the set map
$$\r \colon\spec S \; @> \; P \longmapsto f^{-1}(P) \; >> \; \spec R
$$
is Zariski continuous, and (2)
$$\r^{-1} V_R(X) \; = \; V_S(f(X)).$$

\subhead 1.2 \endsubhead Now let $f\colon R \rightarrow S$ be a
homomorphism of noncommutative rings, and let
$$\r \colon \spec S \longrightarrow \spec R $$
denote the correspondence assigning to each $P \in \spec S$ the set of
prime ideals of $R$ minimal over $f^{-1}(P)$. Adapting \cite{\ArtSch, \S
4}, we will say that $\r$ is {\sl continuous\/} provided ($1'$):
$$\r^{[-1]}V \; \colon= \; \{ P \in \spec S : \r P \subseteq V \}$$
is closed for all closed subsets $V$ of $\spec R$. It need not be true
that $\r$ is continuous, even when $R$ and $S$ are noetherian; see
(2.5). Continuity does hold when $R$ and $S$ satisfy a polynomial
identity; see \cite{\ArtSch, 4.6v} and (2.10).

One generalization of (2) might require that
$$\r^{[-1]} V_R(X) \; = \; V_S(f(X))$$
for all $X \subseteq S$. But it is easy to show that $\r$ can be
continuous while not satisfying this hypothesis; see (2.4iii). Another
possible generalization is ($2'$): For all ideals $I$ of $R$,
$$\r^{[-1]}V_R(I) \; = \; V_S(I^S), \quad \text{where} \quad I^S
\colon= \ann_S \big( S/Sf(I) \big).$$
It follows, for example, from (3.17) that ($2'$) is also strictly
stronger than ($1'$). However, condition ($2'$) will be useful in our
``point free'' approach, described next.

\subhead 1.3 \endsubhead Let $\SPEC R$ denote the category whose
objects are the Zariski closed subsets of $\spec R$ and whose
morphisms are the inclusions; similarly define $\SPEC S$. In \S 5 we
consider the functors
$$\lambda \colon \SPEC S \; @> \; V_S(J) \mapsto V_R(f^{-1}(J)) \; >>
\; \SPEC R, \quad \text{and} \quad \rho \colon \SPEC R \; @> \;
V_R(I) \mapsto V_S(I^S) \; >> \; \SPEC R,$$
where $I$ is a semiprime ideal of $R$ and $J$ is a semiprime ideal of
$S$. When $R$ and $S$ are commutative, it is easy to check that
$\lambda$ is left adjoint to $\rho$; this adjointness amounts,
essentially, to a reformulation of (2).

In our main result, (3.15), we give precise criteria for $\lambda$ to be a
left adjoint to $\rho$, under certain hypotheses (satisfied by left noetherian
rings and affine PI algebras); in particular, this adjointness holds if and
only if $\r$ is a single-valued continuous function (allowing for a slight
abuse of notation) and ($2'$) holds. When $S$ is left noetherian, further
equivalent conditions are given, amounting to a ``nearly centralizing''
property. The moral is that, other than for centralizing extensions, this
adjointness is a rare occurrence.

\subhead 1.4 \endsubhead In the approach to noncommutative algebraic geometry
in \cite{\Rostwo; \VdB}, the ring homomorphism $f\colon R \rightarrow S$
provides only one example of an affine map between affine noncommutative
spaces. Indeed, some of our analysis below can be formulated for more general
morphisms between noncommutative spaces, and a greater portion can be restated
for the setting in which the homomorphism $f\colon R \rightarrow S$ is
replaced by an appropriate $R$-$S$-bimodule. While a few of the definitions
and preliminary results in this paper are presented within this broader
context, we leave a more complete generalization to the interested
reader. Recent studies on noncommutative ring homomorphisms (and
generalizations) from this point of view include \cite{\Smione; \Smitwo;
\SmiZha}.

\subhead 1.5 \endsubhead Our emphasis on categories of closed --
rather than open -- subsets of topological spaces is a matter of
convenience and personal preference. All of the results and
observations below have dual versions involving the categories of open
subsets of topological spaces, with inclusions again providing the
morphisms.

\subhead 1.6 Conventions and Notation \endsubhead (i) Let $A$ be a
ring. We will always assume that the Zariski topology has been applied
to $\spec A$, and we will continue to use the notation $\SPEC$ as in
(1.3). If $I$ is an ideal of $A$, we will use $\sqrt{I}$ to denote the
prime radical of $I$, and if $U$ is a set of prime ideals in $A$ we
will use $I(U)$ to denote their intersection; note that $V_A(I(U))$ is
the closure of $U$ in $\spec A$.

(ii) Let $A$ and $B$ be rings. We will use $_AM$ as an abbreviation
for ``the left $A$-module $M$.'' We will similarly use $M_A$ for right
$A$-modules and $_AM_B$ for $A$-$B$-bimodules.  We will use $\ann _AM$
to denote the annihilator of $_AM$ and $\ann M_A$ to denote the
annihilator of $M_A$. The category of left $A$-modules will be denoted
$\Mod A$.

(iii) The reader is referred to \cite{\GooWar; \McCRob} for further
ring-theoretic background information.

\subhead Acknowledgment \endsubhead I am grateful to F. Van Oystaeyen for
helpful comments on the literature.

\head 2. Continuous Correspondences \endhead

In this section we consider ring homomorphisms and continuity. This
discussion can be regarded as a continuation of \cite{\ArtSch, \S 4},
where correspondences between the maximal spectra of affine PI
algebras are considered.

Throughout this section, $f\colon R \rightarrow S$ will be a
homomorphism of rings. 

\subhead 2.1 \endsubhead (i) Let $X$ and $Y$ be sets. By a {\sl
correspondence\/} $\c\colon X \rightarrow Y$ we mean a function from
$X$ into the set of subsets of $Y$. Following common practice, we will
define
$$\c U \colon= \bigcup_{u \in U}\c u \qquad \text{and} \qquad
\c^{-1} V \colon= \{ u \in U : \c u \cap V \ne \emptyset \},$$
for subsets $U$ of $X$ and $V$ of $Y$. However, it will be more
convenient for our purposes to use the following generalization of the
inverse of a function,
$$\c^{[-1]} V \colon= \{ u \in U : \c u \subseteq V \}. $$
Note that $X - \c^{-1}V = \c^{[-1]}(Y-V)$.  Also,
$$\c^{[-1]}\c U \supseteq U, \quad \c\c^{[-1]}V \subseteq V,
\quad \c U \subseteq \c U', \quad \text{and} \quad \c^{[-1]} V
\subseteq \c^{[-1]} V' ,$$
for all $U \subseteq U' \subseteq X$ and $V \subseteq V' \subseteq Y$.

(ii) Let $X$ and $Y$ be topological spaces. Following \cite{\ArtSch,
\S 4}, we will say that the correspondence $\c\colon X \rightarrow Y$
is {\sl continuous\/} provided $\c^{-1}W$ is open for all open subsets
$W$ of $Y$, or equivalently, provided $\c^{[-1]}Z$ is closed for all
closed subsets $Z$ of $Y$.

\subhead 2.2 \endsubhead The correspondences of spectra of interest to
us appear within the following more general framework. Let $\alpha \colon
\Mod B \rightarrow \Mod A$ be a covariant functor, for
rings $A$ and $B$. Given an ideal $J$ of $B$, set
$$J^\alpha \colon= \ann_A \alpha (B/J).$$
We obtain a correspondence $\r(\alpha) \colon \spec B \rightarrow \spec
A$, sending each $P \in \spec B$ to the set of prime ideals
of $A$ minimal over $P^\alpha$. (It may be the case that
$P^\alpha = A$, in which case $\r(\alpha) P$ will be empty. However,
using Zorn's lemma, if $J$ is an ideal of $A$ contained within at
least one $Q \in \spec A$, then there exists a $Q' \in \spec
A$ such that $Q' \subseteq Q$ and such that $Q'$ is minimal over
$J$.)

\subhead 2.3 \endsubhead Applying (2.2) to the restriction of scalars
functor $\Mod S \rightarrow \Mod R$, we obtain the correspondence
(which we will denote) $\r \colon \spec S \rightarrow \spec R$,
sending each $P \in \spec S$ to the nonempty set
$$\{ Q \in \spec R : \text{$Q$ is minimal over $f^{-1}(P) = \ann
_R(S/P)$} \}.$$
If $Q$ is a prime ideal of $R$, then $\r^{-1}Q$ is commonly referred
to as the set of prime ideals of $S$ ``lying over'' $Q$.

\subhead 2.4 \endsubhead Let $I$ be an ideal of $R$.

(i) Note that
$$\r ^{[-1]} V_R(I) = \{ P \in \spec S : \sqrt{f^{-1}(P)} \supseteq I
\}.$$

(ii) When $R$ and $S$ are commutative, $\r$ is the continuous
function from $\spec S$ to $\spec R$ mapping each prime ideal $P$ of
$S$ to the prime ideal $f^{-1}(P)$ of $R$, and
$$\r ^{[-1]}V_R(I) = \r^{-1}V_R(I) = V_S(f(I)).$$

(iii) When $S$ is not commutative, the equality in (ii) need
not hold. For example, set
$$S = \bmatrix k & k \\ k & k \endbmatrix, \quad R = \left\{ \bmatrix
\alpha & \beta \\ 0 & \alpha \endbmatrix : \alpha, \beta \in k
\right\} \subseteq S, \quad \text{and} \quad I = \bmatrix 0 & k \\ 0 &
0 \endbmatrix \subseteq R,$$
where $k$ is a field. Let $f$ be the inclusion of $R$ in $S$. Then
$\spec S = \{ 0 \}$, $\spec R = \{ I \}$, $\r$ is continuous, and
$$ \r^{[-1]}V_R(I) = \{ 0 \} \ne \emptyset = V_S(f(I)) .$$

(iv) When $\sqrt{f^{-1}(P)}$ is nilpotent modulo $f^{-1}(P)$ for all
$P \in \spec S$, the equality in (ii) can be replaced by 
$$\align \r ^{[-1]} V_R(I) &= \{ P \in \spec S : \text{$f^{-1}(P)
\supseteq I^t$ for some positive integer $t$} \} \\ &= \{ P \in \spec
S : \text{$P \supseteq f(I)^t$ for some positive integer $t$} \} \\ &=
\bigcup _{t \geq 1} V_S(f(I)^t) . \endalign$$

\subhead 2.5 \endsubhead We can see as follows that $\r$ need not be
continuous, even when $R$ and $S$ are noetherian.

Let $k$ be a field of characteristic zero, and suppose that $S$ has
been chosen to be the enveloping algebra of $\frak{sl}_2(k)$. Let $\{
E, F, H \}$ be the standard $k$-basis for $\frak{sl}_2(k) \subset S$
(cf., e.g., \cite{\Dix, \S 1.8}), with $[H,E] = 2E$, $[H,F] = -2F$, and
$[E,F] = H$. Assume that $R = k\{ E\} \subset S$ and that $f$ is the
inclusion map. Let $I = \langle E \rangle$. It is well known that $R$
is a polynomial ring in $E$ and that $S$ is noetherian.  Moreover, if
$P$ is the kernel of a finite dimensional irreducible representation
of $S$, then $I^t \subset P$ for some positive integer $t$. (This last
assertion immediately follows, e.g., from \cite{\Dix, \S 1.8}.)

We can now see that $U = \r ^{[-1]}V_R(I) \subset \spec S$ contains the kernel
of every finite dimensional irreducible representation of $S$. It is well
known that the intersection of these kernels is zero. Therefore, $I(U) = 0$, a
prime ideal of $S$. However, the ideal $0$ of $S$ cannot be contained in $U$,
and so $U \ne V_S(I(U))$. Therefore, $U$ is not closed, and $f$ is not
continuous.

\subhead 2.6 \endsubhead Continuity does hold in the following
commonly occurring special case: Suppose that $f^{-1}(P)$ is a
semiprime ideal of $R$ for every prime ideal $P$ of $S$. (See, e.g.,
\cite{\McCRob, Chapter 10} for settings in which this hypothesis
holds.) Then, if $I$ is an ideal of $R$,
$$\r ^{[-1]}V_R(I) = \{ P \in \spec S : f^{-1}(P) \supseteq I \} = \{
P \in \spec S : P \supseteq f(I)\} = V_S(f(I)).$$
Hence $\r$ is continuous.

\subhead 2.7 \endsubhead In the remainder of this section we establish
continuity in the presence of a bound on Goldie ranks.

(i) Let $A$ be a ring for which every prime factor is left or right
Goldie. Set
$$\spec _n A = \{ P \in \spec A : \rank (A/P) \leq n \},$$
where ``rank'' means ``Goldie rank'' and where $n$ is a positive
integer. Equip $\spec _n A$ with the relative Zariski topology.

(ii) Suppose that all of the prime factors of $R$ and $S$ are left or right
Goldie. Let $P \in \spec_n S$. It follows from \cite{\War} that $\r P \in
\spec _n R$.

\proclaim{2.8 Lemma} Let $A$ be a subring of a prime left or right
Goldie ring $B$. Suppose that the Goldie rank of $B$ is $t$, and let
$N$ denote the prime radical of $A$. Then $N^t = 0$. \endproclaim

\demo{Proof} Let $F$ be the Goldie quotient ring of $B$. By
assumption, $F$ has length $t$ as a left $F$-module, and so there
exists an $F$-$A$-bimodule composition series
$$0 = F_0 \subset F_1 \subset \cdots \subset F_s = F,$$
for some $s \leq t$. For $1 \leq i \leq s$, set
$$Q_i = \ann (F_i/F_{i-1})_A .$$
Then $F.Q_s\cdots Q_1 = 0$, and it is easy to check that
$Q_1,\ldots,Q_s$ are prime ideals of $A$. In particular, $Q_s\cdots
Q_1 = 0$ in $A$, and so $N^t \subseteq N^s = 0$. \qed\enddemo

\proclaim{2.9 Proposition} Let $n$ be a positive integer, and assume
that all of the prime factors of $R$ and $S$ are left or right
Goldie. Then $\r \colon \spec _n S \rightarrow \spec _n R$ is
continuous.
\endproclaim

\demo{Proof} Without loss of generality, we may assume that $R$ is a
subring of $S$ and that $f$ is the inclusion map. Let $I$ be an ideal
of $R$, and set $V = V_R(I)$. It now follows from (2.8), and our
earlier observations, that
$$\multline \left(\r ^{[-1]}\left(V \cap \spec_n R\right)\right) \cap
\spec_n S = \left(\r ^{[-1]}V\right) \cap \spec_n S = \\ \{ P \in
\spec_n S : P \supseteq I^n \} = \big( V_S(
I^n) \big) \cap \spec_n S. \endmultline $$
The proposition follows. \qed\enddemo

\proclaim{2.10 Corollary} {\rm (cf\. \cite{\ArtSch, 4.6v})} If
$S$ is a PI ring then $\r$ is continuous. \endproclaim

\demo{Proof} Assume that $S$ is PI. It follows from Posner's theorem
that every prime factor of $R$ and $S$ is Goldie. It follows from
basic PI theory that there exists a finite upper bound for the Goldie
ranks of the prime factors of $S$. The corollary now follows from
(2.9). \qed\enddemo

\subhead 2.11 \endsubhead In \cite{\ArtSch, 4.6v} it is noted that the
correspondence $\r\colon \max S \rightarrow \max R$ is continuous when
$R$ and $S$ are PI algebras affine over a field. However, the proof
given there (in the last paragraph on page 307) appears to be
incorrect.

\subhead 2.12 \endsubhead In \cite{\ArtSch, 4.7} it is shown that the
homomorphism $f\colon R \rightarrow S$ can be chosen with the
following properties: (i) $R$ and $S$ are PI algebras affine over a
field, (ii) there exists a closed subset $V$ of $\spec R$ for which
$\r^{-1}V$ is not closed in $\spec S$. As noted in \cite{\ArtSch, 4.7},
it follows that ``$\r^{-1}(\text{open})$ is open'' continuity does not
imply ``$\r^{-1}(\text{closed})$ is closed'' continuity.

\subhead 2.13 \endsubhead We ask: (i) Must $\r$ be continuous when $S$
is FBN? (ii) Must $\r$ be continuous when $S$ is finitely generated as
an $R$-module?

\head 3. Adjointness \endhead

Throughout this section, $f\colon R \rightarrow S$ will be a ring
homomorphism, and $\r$ will denote the correspondence from $\spec S$ to $\spec
R$ described in (2.3). In our main result, (3.15), we determine -- under
addtional hypotheses introduced in (3.7) -- when adjointness holds for the
functors, between $\SPEC R$ and $\SPEC S$, arising from restriction and
extension of scalars.

We begin with some preliminaries on functors, correspondences, and
topological spaces.

\subhead 3.1 \endsubhead Let $X$ be a topological space, and let
$\closed X$ denote the category whose objects are the closed subsets
of $X$ and whose morphisms are the inclusions. If $U$ is a subset of
$X$, we will denote the closure of $U$ in $X$ by $\overline{U}$.

\subhead 3.2 \endsubhead Let $X$ and $Y$ be topological spaces.

(i) Let $\varphi $ be a covariant functor from $\closed X$ to $\closed
Y$, and let $\psi$ be a covariant functor from $\closed Y$ to $\closed
X$. Then $\varphi $ is a left adjoint to $\psi$ exactly when
$$\varphi U \subseteq V \Longleftrightarrow U \subseteq \psi V,$$
for all $U \in \closed X$ and $V \in \closed Y$. Now suppose that
$\psi$ and $\psi'$ are both right adjoints to $\varphi $, and let $V \in
\closed Y$. Then
$$\psi V \subseteq \psi V \Rightarrow \varphi \psi V \subseteq V
\Rightarrow \psi V \subseteq \psi' V .$$
Similarly, $\psi' V \subseteq \psi V$. It follows that $\psi$ and
$\psi'$ must be the same functor.

(ii) Let $\c\colon X \rightarrow Y$ be a (not necessarily
continuous) correspondence. We obtain covariant functors
$$\varphi^{\c} \colon\closed X @> U \mapsto \overline{\c U} >> \closed
Y, \quad \text{and} \quad \varphi_{\c} \colon\closed Y @> V \mapsto
\overline{\c^{[-1]}V} >> \closed X.$$
Moreover, $\varphi^{\c}$ is a left adjoint to $\varphi_{\c}$ exactly when
$$\c U \subseteq V \Leftrightarrow U \subseteq
\overline{\c^{[-1]}V} $$
for all closed subsets $U$ of $X$ and $V$ of $Y$. Consequently, if
$\c$ is continuous, it immediately follows that $\varphi^{\c}$ is a
left adjoint to $\varphi_{\c}$. Conversely, if $\varphi^{\c}$ is a left
adjoint to $\varphi_{\c}$, then
$$\overline{\c^{[-1]}V} \subseteq \overline{\c^{[-1]}V}
\Rightarrow \c\left(\overline{\c^{[-1]}V}\right) \subseteq V
\Rightarrow \overline{\c^{[-1]}V} \subseteq
\c^{[-1]}\c\left(\overline{\c^{[-1]}V}\right) \subseteq
\c^{[-1]}V.$$
We conclude that $\varphi^{\c}$ is a left adjoint to $\varphi_{\c}$ if and
only if $\c$ is continuous. 

\subhead 3.3 \endsubhead We now introduce functors between spectra in
a somewhat more general framework. Assume that $A$ and $B$
are rings, and that $\alpha \colon \Mod B \rightarrow \Mod A$
is a covariant functor. Recall the notation of (2.2).

(ii) Following (3.2ii), we obtain the functors
$$\varphi^{\r(\alpha)} \colon \SPEC B \rightarrow \SPEC A
\quad \text{and} \quad \varphi_{\r(\alpha)} \colon \SPEC A
\rightarrow \SPEC B .$$

(ii) Suppose that $\alpha$ is right exact. Then the assignment $J
\mapsto J^\alpha$ preserves inclusions, and thus induces a functor
$$\theta^\alpha \colon \SPEC B \; @> \; V \mapsto V_A\big(
I(V)^\alpha \big) \; >> \; \SPEC A .$$

\subhead 3.4 \endsubhead Retain the notation of (3.3), and assume that there
exists an $A$-$B$-bimodule $M$ such that $\alpha L = M\otimes_B L$, for each
left $B$-module $L$. Recall, by Watts' theorem (see, for example, \cite{\Ste,
IV.10.1}, that this assumption holds if and only if $\alpha$ posesses a right
adjoint.

(i) Observe that
$$J^\alpha = \ann _A \big( M/M.J \big),$$
for all ideals $J$ of $B$.

(ii) Note, for ideals $J_1$ and $J_2$ of $B$, that
$$J_1^\alpha J_2^\alpha.M \subseteq J_1^\alpha.M.J_2 \subseteq
M.J_1J_2 ,$$
and so $J_1^\alpha J_2^\alpha \subseteq (J_1J_2)^\alpha$.

(iii) Let $J$ be an ideal of $B$, and suppose that $Q$ is a prime
ideal of $A$ containing $J^\alpha$. Using Zorn's lemma, we can choose
an ideal $P$ of $B$ maximal such that $P \supseteq J$ and such that
$P^\alpha \subseteq Q$; it follows from (ii) that $P$ must be
prime. Therefore,
$$Q \supseteq P^\alpha \supseteq \left(\sqrt{J}\right)^\alpha, \quad
\text{and so} \quad \sqrt{J^\alpha} \supseteq
\left(\sqrt{J}\right)^\alpha \supseteq J^\alpha.$$
It follows that
$$\theta^\alpha V_B(J) = V_A\big (J^\alpha),$$
for all ideals $J$ of $B$. 

(iv) Let $J$ be an ideal of $B$, and set
$$X = \r(\alpha) V_B(J) \subseteq \spec A.$$
If $Q \in X$ then $Q \supseteq J^\alpha$, and so $\overline{X}
\subseteq V_A(J^\alpha)$. Conversely, choose $Q \in V_A(J^\alpha)$. As
in (iii), there exists a prime ideal $P$ of $B$ such that $P \supseteq
J$ and such that $Q \supseteq P^\alpha$. There then exists (by another
Zorn's lemma argument) a prime ideal $Q'$ of $A$ minimal over
$P^\alpha$ such that $Q' \subseteq Q$. Because $Q' \in X$, we see that
$Q \in \overline{X}$, and so
$$\varphi^{\r(\alpha)}V_B(J) = \overline{X} = V_A(J^\alpha) =
\theta^\alpha V_B(J) .$$
We see, in the present setting, that $\theta^\alpha$ and
$\varphi^{\r(\alpha)}$ are the same functor.

\subhead 3.5 \endsubhead Applying (3.4) to the restriction of scalars
functor $\Mod S \rightarrow \Mod R$, we obtain the functor $\lambda
\colon \SPEC S \rightarrow \SPEC R$, sending
$$V_S(J) \longmapsto V_R(f^{-1}(J)),$$
for ideals $J$ of $S$. Again using (3.4), we see that $\lambda =
\varphi^{\r}$.

\subhead 3.6 \endsubhead (i) For each ideal $I$ of $R$, set
$$I^S = \ann_S \big(S/Sf(I) \big) .$$
Applying (3.4) to the extension of scalars functor $\Mod R \rightarrow
\Mod S$, we now obtain the functor $\rho \colon \SPEC R
\rightarrow \SPEC S$, sending
$$V_R(I) \longmapsto V_S\big( I^S \big),$$
for ideals $I$ of $R$.

(iii) Suppose that $R$ and $S$ are commutative. Then $\r\colon\spec S
\rightarrow \spec R$ is a continuous function, and, in the notation of
(3.4), $\rho = \varphi_{\r}$. Moreover, following (3.2ii) we see that
$\lambda$ is a left adjoint to $\rho$. 

\subhead 3.7 \endsubhead For the remainder of this section we will
assume that (i) all semiprime factors of $R$ and $S$ are left or right
Goldie, and (ii) the prime radicals of all of the factors of $R$ and $S$
are nilpotent.

\subhead 3.8 \endsubhead (i) The hypotheses in (3.7) hold, of course,
when $R$ and $S$ are left or right noetherian. 

(ii) Suppose that $R$ and $S$ are each affine over a commutative
noetherian ring and satisfy a polynomial identity. Then (3.7i) follows
from Posner's theorem, and (3.7ii) follows from \cite{\Bra}.

(iii) Let $I$ be an ideal of $R$ or $S$. It follows from (3.7i) that
$\sqrt{I}$ is the intersection of finitely many prime ideals and then
from (3.7ii) that $I$ contains a finite product of prime ideals. In
particular, there are finitely many prime ideals minimal over $I$.

\subhead 3.9 \endsubhead In (3.10) through (3.14) we will further assume
that $R$ is a subring of $S$ and that $f$ is inclusion. 

\proclaim{3.10 Lemma} {\rm (i)} If $J$ is an ideal of $S$ then
$\lambda V_S(J) = V_R(J\cap R)$. {\rm (ii)} If $I$ is an ideal of $R$ then
$\rho V_R(I) = V_S(I^S)$. \endproclaim

\demo{Proof} (i) Let $J$ be an ideal of $S$. For sufficiently large
$t$,
$$\left(\sqrt{J}\cap R\right)^t \subseteq J\cap R \subseteq \sqrt{J}
\cap R,$$
and so
$$V_R(J\cap R) = V_R\left(\sqrt{J} \cap R\right) = \lambda V_S(J) .$$

{\rm (ii)} Let $I$ be an ideal of $R$. By (3.4ii), for sufficiently
large $t$,
$$\left(\left(\sqrt{I }\right)^S\right)^t \subseteq
\left(\left(\sqrt{I}\right)^t\right)^S \subseteq I^S \subseteq
\left(\sqrt{I}\right)^S,$$
and so
$$V_S\left(I^S\right) = V_S\left(\left(\sqrt{I}\right)^S\right) =
\rho V_R(I) . $$
(The preceding two arguments are symmetrical -- note that $J \cap R =
\ann_R (S/SJ)$, for ideals $J$ of $S$.)  \qed\enddemo

\subhead 3.11 \endsubhead We can now see, in the present situation,
that $\lambda$ is a left adjoint to $\rho$ exactly when
$$V_S(J) \subseteq V_S(I^S) \quad \Longleftrightarrow \quad V_R(J\cap
R) \subseteq V_R(I), $$
or equivalently,
$$I^S \subseteq \sqrt{J} \quad \Longleftrightarrow \quad I \subseteq
\sqrt{J\cap R},$$
for all ideals $I$ of $R$ and $J$ of $S$.

\proclaim{3.12 Lemma} {\rm (i)} Let $I$ be an ideal of $R$, let $J$ be
an ideal of $S$, and suppose that $V_R(J\cap R) \subseteq V_R(I)$.
Then $V_S(J) \subseteq V_S(I^S)$. 

{\rm (ii)} $\lambda$ is a left adjoint to $\rho$ if and only if
$$V_S(J) \subseteq V_S(I^S) \quad \Longrightarrow \quad V_R(J\cap R)
\subseteq V_R(I),$$
for all ideals $I$ of $R$ and $J$ of $S$. \endproclaim

\demo{Proof} (i) Since $I \subseteq \sqrt{J\cap R}$, there exists a
positive integer $t$ such that $I^t \subseteq J\cap R$. Hence $I^tS
\subseteq J$, and so $(I^t)^S \subseteq J$. Therefore, by (3.10),
$$V_S(J) \subseteq V_S\left(\left(I^t\right)^S\right) =
\rho V_R\left(I^t\right) = \rho V_R(I) = V_S(I^S). $$

(ii) Follows immediately from (i) and (3.11). \qed\enddemo

\proclaim{3.13 Lemma} The following are equivalent.

{\rm (i)} $\lambda$ is a left adjoint to $\rho$.

{\rm (ii)} For all $P \in \spec S$ and $Q \in \spec R$,
$$Q^S \subseteq P \quad \Longrightarrow \quad Q \subseteq
\sqrt{P\cap R} .$$
\endproclaim

\demo{Proof} It follows immediately from (3.11) that (i) $\Rightarrow$
(ii).

Conversely, assume that (ii) is true, that $I$ is an ideal of $R$,
that $J$ is an ideal of $S$, and that $V_S(J) \subseteq
V_S(I^S)$. Then $I^S \subseteq \sqrt{J}$. Let $P$ be a prime ideal of
$S$ minimal over $J$.

Using Zorn's lemma we can choose an ideal $Q$ of $R$ maximal
among the ideals $I'$ of $R$ for which $I' \supseteq I$ and $I^S
\subseteq P$. Because $P$ is prime, (3.4ii) ensures that $Q$ is
prime. Therefore, by assumption, $Q \subseteq \sqrt{P\cap R}$, and so
$I \subseteq \sqrt{P\cap R}$. Consequently, $I^t \subseteq P\cap R$
for a sufficiently large positive integer $t$.

Since $P$ was arbitrarily chosen among the finitely many prime ideals
of $S$ minimal over $J$, we see that $I^t \subseteq \sqrt{J}\cap R$
for sufficiently large $t$. However, $(\sqrt{J}\cap R)^t \subseteq
J\cap R$ for sufficiently large $t$, and so $I^t \subseteq J\cap R$
for sufficiently large $t$. Therefore, $I \subseteq \sqrt{J\cap R}$.
Hence $V_R(J\cap R) \subseteq V_R(I)$, and it follows from (3.12ii)
that (ii) $\Rightarrow$ (i). \qed\enddemo

\proclaim{3.14 Lemma} Let $P \in \spec S$. Then there exists a $Q \in
\spec R$ such that $Q$ is minimal over $P\cap R$ and such that $Q^S
\subseteq P$. \endproclaim

\demo{Proof} We may assume, without loss of generality, that $P = 0$.
Next, by (3.8iii), there exists a prime ideal $\widehat{Q}$ of $R$
such that $\widehat{Q}.N = 0$ for some nonzero ideal $N$ of
$R$. Choose a minimal prime ideal $Q$ of $R$ such that $Q \subseteq
\widehat{Q}$, and let $F$ denote the Goldie quotient ring of
$S$. Since $F.Q.N = 0$, and since $\ann F_S = 0$, we see that $F.Q
\ne F$. Consequently, $F/FQ$ is a nonzero $F$-$R$-bimodule. By
Goldie's theorem, every left $S$-submodule of $F/FQ$ must have
annihilator equal to $P$.

Now note that $F/FQ$ contains a nonzero $S$-$R$-bimodule factor of
$S/SQ$. In particular, there exists an $S$-$R$-bimodule factor $B$ of
$S/SQ$ with $\ann_SB = 0$. Thus $Q^S = \ann_S(S/SQ) = 0$, and the
lemma follows. \qed\enddemo

\proclaim{3.15 Theorem} Assume that $f\colon R \rightarrow S$ be a ring
homomorphism, that all semiprime factors of $R$ and $S$ are left or right
Goldie, and that the prime radicals of all of the factors of $R$ and $S$ are
nilpotent.

{\rm (1)} The following are equivalent.

{\rm (i)} $\lambda$ is a left adjoint to $\rho$.

{\rm (ii)} The canonical correspondence $\r\colon\spec S \rightarrow
\spec R$ defined in {\rm (2.3)\/} is a single-valued continuous function, and
$$\r ^{[-1]}V_R(I) = V_S(I^S),$$
for all ideals $I$ of $R$. 

\noindent {\rm (2)} If $S$ is left noetherian then {\rm (i)}, {\rm
(ii)} and the following are equivalent.

{\rm (iii)} For each $Q \in \spec R$ there is a positive integer
$t$ such that $f(Q)^tS \subseteq Sf(Q)$. 

{\rm (iv)} For each ideal $I$ of $R$ there is a positive integer $t$
such that $f(I)^tS \subseteq Sf(I)$. 
\endproclaim

\demo{Proof} We may assume, without loss of generality, that $R$ is a
subring of $S$ and that $f$ is inclusion.

(1) (i) $\Rightarrow$ (ii): Let $P \in \spec S$. By (3.14), we can
choose $Q \in \spec R$ such that $Q$ is minimal over $P\cap R$ and
such that $Q^S \subseteq P$. By (3.11), $Q \subseteq \sqrt{P\cap R}$,
and so $Q = \sqrt{P\cap R}$. Hence $\r P = \{ Q \}$, and $\r$ is a
single-valued function.

Now let $I$ be an ideal of $R$, and note that $P \in \r ^{[-1]}V_R(I)$
if and only if $I \subseteq \sqrt{P\cap R}$. Hence, by (3.11),
$\r ^{[-1]}V_R(I) = V_S(I^S)$. In particular, $\r$ is continuous.

(ii) $\Rightarrow$ (i): Assume that $P \in \spec S$, that $Q \in \spec R$, and
that $Q^S \subseteq P$. In other words, $P \in V_S(Q^S)$. By hypothesis,
$V_S(Q^S) = \r ^{[-1]}V_R(Q)$, and hence $P \in \r ^{[-1]}V_R(Q)$. Therefore,
$\r P \subseteq V_R(Q)$, and so $Q \subseteq \sqrt{P\cap R}$. It now follows
from (3.13) that $\lambda$ is a left adjoint to $\rho$.

(2) Assume that $S$ is left noetherian.

(i) $\Rightarrow$ (iii): Suppose that $S/SQ \ne 0$; the desired
conclusion immediately holds true otherwise. Next, since $S$ is left
noetherian, there exists a series of $S$-$R$-bimodules,
$$0 = M_0 \subset M_1 \subset \cdots \subset M_n = S/SQ,$$
such that for each $1 \leq i \leq n$,
$$P_i = \ann_S(M_i/M_{i-1}) \in \spec S$$
(see, e.g., \cite{\GooWar, 2.13}). In particular, $Q^S$ is contained in
each of $P_1,\ldots,P_n$. In view of (3.11), it now follows from our
assumptions that $Q \subseteq \sqrt{P_i\cap R}$, for $1 \leq i \leq
n$. Therefore, for sufficiently large $t$, $Q^t \subseteq P_1\cdots
P_n$. Consequently, $Q^t.(S/SQ) = 0$, and so $Q^tS \subseteq SQ$.

(iii) $\Rightarrow$ (i): Assume that $Q \in \spec R$, that $P \in \spec
S$, and that $Q^S \subseteq P$. Choose $t$ such that $Q^tS \subseteq
SQ$. Then $SQ^tS \subseteq SQ$, and so $SQ^tS \subseteq
Q^S$. Hence, $Q^t \subseteq (SQ^tS)\cap R \subseteq P\cap
R$. Therefore, $Q \subseteq \sqrt{P\cap R}$. By (3.13), $\lambda$ is a left
adjoint to $\rho$. 

(iii) $\Leftrightarrow$ (iv): Assume (iii), and let $I$ be an arbitary ideal
of $R$. Choose $Q_1,\ldots, Q_n \in \spec R$ such that $\sqrt{I} = Q_1 \cap
\cdots \cap Q_n$ and such that $Q_1\ldots Q_n \subseteq I$. Then, by
assumption, for a sufficiently large positive integer $t$, $I^ntS \subseteq
SQ_1\cdots Q_n \subseteq SI$, and (iv) holds true. The converse is
trivial. \qed\enddemo

\subhead 3.16 \endsubhead It is easy to see that the conditions (iii) and (iv)
of (3.15) are satisfied when the homomorphism $f\colon R \rightarrow S$ is
centralizing (i.e., $S$ is generated as a left $R$-module by a set $X$ such
that $r.x = x.r$ for all $r \in R$ and $x \in X$). Non-centralizing
homomorphisms for which (3.15iii, iv) hold are more rare, although ring
embeddings associated to nilpotent Lie superalgebras provide such examples;
see \cite{\Letone; \Lettwo} for details. We can view ring homomorphisms
satisfying (3.15iv) as being ``nearly centralizing.''

\subhead 3.17 \endsubhead It is not true that $\lambda$ is a left
adjoint to $\rho$ if and only if $\r$ is a single-valued continuous
function. To provide an easy illustration, let $k$ be a field of
characteristic zero and let $S$ denote the first Weyl algebra over
$k$: $S$ is generated by $x$ and $y$, subject only to the relation $yx
- xy = 1$.  Let $R$ be the commutative polynomial ring $k[x]$,
identified with the subalgebra of $S$ generated by $x$, and let $f$
denote the inclusion homomorphism.

Let $P$ denote the zero ideal of $S$. Then $\spec S = \{ P \}$ and
$P\cap R \in \spec R$. Hence $\r$ is a single-valued continuous
function.

Now let $I$ be the ideal of $R$ generated by $x$. Then $SI = Sx$ is a
proper left ideal of $S$, and so $S/SI \ne 0$. Since $S$ is a simple
ring, $I^S = 0$. Also, $I^S \subset P$ and $I \not\subseteq
\sqrt{P\cap R}$. Therefore, by (3.11), $\lambda$ is not a left adjoint
to $\rho$.

\enddocument